\documentclass[11pt]{article}
\usepackage{amscd}
\usepackage{amsfonts}
\usepackage{amsmath}
\usepackage{amssymb}
\usepackage{amsthm}
\usepackage{bbm}
\usepackage{CJK}
\usepackage{fancyhdr}
\usepackage{graphicx}
\usepackage{hyperref}
\usepackage{indentfirst}
\usepackage{latexsym}
\usepackage{mathrsfs}
\usepackage{xypic}

\usepackage[top=1in,bottom=1in,left=1.25in,right=1.25in]{geometry}
\def\<{\langle}
\def\>{\rangle}

\def\c{\cdot}

\def\o{\otimes}

\date{}
\begin{document}
\renewcommand{\baselinestretch}{1.2}
\renewcommand{\arraystretch}{1.0}
\title{\bf Derivations and deformations of $\delta$-Jordan Lie supertriple systems}
\author{{\bf Shengxiang Wang$^{1}$, Xiaohui Zhang$^{2}$, Shuangjian Guo$^{3}$\footnote
        {Corresponding author(Shuangjian Guo): shuangjianguo@126.com} }\\
1.~ School of Mathematics and Finance, Chuzhou University,\\
 Chuzhou 239000,  China \\
2.~  School of Mathematical Sciences, Qufu Normal University, \\Qufu 273165, China.\\
 3.~ School of Mathematics and Statistics, Guizhou University of\\ Finance and Economics,Guiyang 550025,China.}
 \maketitle
\begin{center}
\begin{minipage}{13.cm}

{\bf \begin{center} ABSTRACT \end{center}}
Let $T$  be a $\delta$-Jordan Lie supertriple system.
We first introduce the notions of generalized derivations and representations of $T$ and present  some properties.
Also,  we study the low dimension cohomology and the coboundary operator of $T$, and then we  investigate the deformations and  Nijenhuis operators of $T$ by choosing some suitable cohomology.
 \smallskip

 \smallskip

 \end{minipage}
 \end{center}
 \normalsize\vskip0.5cm

\section{Introduction}
\def\theequation{\arabic{section}. \arabic{equation}}
\setcounter{equation} {0}

Lie triple systems arosed initially in Cartan's study of Riemannian geometry.
Jacobson \cite{Jacobson} first introduced Lie triple systems and Jordan triple systems in connection with problems from Jordan theory
and quantum mechanics, viewing Lie triple systems as subspaces of Lie algebras that are closed relative to the ternary product.
Lister \cite{Lister} investigated notions of the radical, semi-simplicity and solvability as
defined for Lie triple systems, and determined all simple Lie triple systems over an algebraically closed field.
Later, the  representation theory,  the central extension, the deformation theory,   bilinear forms and the generalized derivation
of Lie triple systems and Jordan triple systems
have been developed, see \cite{Bremner, Harris, Hodge, Kubo, Stitzinger, Yamaguti, Zhang2009, Zhang2014, Zhang2002}
 \medskip

In \cite{Okubo1997}, Okubo and Kamiya introduced the notion of $\delta$-Jordan Lie triple system,  where $\delta=\pm 1$,
 which is a generalization of both Lie triple systems ($\delta= 1$) and Jordan Lie triple systems ($\delta=-1$).
Later, Kamiya and Okubo \cite{Okubo2002} studied a construction of simple Jordan superalgebras from certain
triple systems.
Recently, Ma and Chen \cite{chen2017} discussed the cohomology theory, the deformations,
Nijenhuis operators, abelian extensions and T$^\ast$-extensions of $\delta$-Jordan Lie triple system.

As a natural generalization of Lie triple systems, Okubo  \cite{Okubo1994} introduced the notion of Lie supertriple systems
 in the study of Yang-Baxter equations.
  Lie supertriple systems have many applications  in high energy physics, and
many important results on Lie supertriple systems have been obtained, see \cite{chen2013, Okubo1994, Okubo2002, chen2015}.
In \cite{Okubo1997}, Okubo and Kamiya introduced the notion of $\delta$-Jordan Lie supertriple system (they still call it Jordan Lie triple system), they presented some nontrivial examples and discussed their quasiclassical property.
In the present paper, we hope to study  generalized derivations,  cohomology theories and   deformations of $\delta$-Jordan Lie supertriple systems.
 \medskip

This paper is organized as follows.
In Section 2, we recall the definition of $\delta$-Jordan Lie supertriple systems and
 construct a kind of  $\delta$-Jordan Lie supertriple systems.
Also, we study generalized derivation algebra of a $\delta$-Jordan Lie supertriple system.
In Section 3, we introduce  notions of the representation and low dimension cohomology of  a $\delta$-Jordan Lie supertriple system.
In Section 4, we consider  the theory of deformations of a $\delta$-Jordan Lie supertriple system  by choosing a suitable cohomology.
In Section 5,  we study  Nijenhuis operators for a $\delta$-Jordan Lie supertriple system  to describe trivial deformations.

\section{Generalized derivations of $\delta$-Jordan Lie supertriple systems}
\def\theequation{\arabic{section}. \arabic{equation}}
\setcounter{equation} {0}

In this section, we start by recalling the definition of $\delta$-Jordan Lie supertriple systems,
then we study its generalized derivations.
\medskip

\noindent{\bf Definition 2.1.} (\cite{Okubo1997})
A $\delta$-Jordan Lie supertriple system is a  $Z_2$-graded vector space $T$
 together with a  triple linear product $[\c,\c,\c]:T\o T\o T\rightarrow T$ satisfying
\begin{eqnarray}
&&(1)~|[a,b,c]|=(|a|+|b|+|c|)(mod~ 2);\\
&&(2)~[b,a,c]=-\delta(-1)^{|a||b|}[a,b,c];\\
&&(3)~(-1)^{|a||c|}[a,b,c]+(-1)^{|b||a|}[b,c,a]+(-1)^{|c||b|}[c,a,b]=0;\\
&&(4)~[a,b,[c,d,e]]=[[a,b,c],d,e]+(-1)^{|c|(|a|+|b|)}[c,[a,b,d],e]\nonumber\\
&&~~~~~~~~~~~~~~~~~~~~~~~~~~~~~~~~~~~~~~~~+\delta(-1)^{(|a|+|b|)(|c|+|d|)}[c,d,[a,b,e]],
\end{eqnarray}
for all $a,b,c,d,e\in T$, where $\delta=\pm 1$ and $|a|$ denotes the degree of the element $a\in T$.
\medskip

\noindent{\bf Remark 2.2.}
Clearly, $T_{\overline{0}}$ is an ordinary $\delta$-Jordan Lie triple system in \cite{chen2017}.
Especially, the case of $\delta= 1$ defines a Lie supertriple system while the other case of $\delta=-1$
may be termed an anti Lie supertriple system as in \cite{Kamiya1988}.
\smallskip

\noindent{\bf Example 2.3.} (\cite{Okubo1997})
Let $(T,[\c,\c])$  be a $\delta$-Jordan Lie superalgebra.
Then $(T,[\c,\c,\c])$  becomes a $\delta$-Jordan Lie supertriple system, where
$[a,b,c]=[[a,b],c]$, for all  $a,b,c\in T$.

\noindent{\bf Example 2.4.}
Let $T$  be a $\delta$-Jordan Lie supertriple system and $t$ an indeterminate.
Set $T'=\{\sum_{i\geq 0}x\o t^i|x\in T\}$, then  $T'$ is a  $\delta$-Jordan Lie supertriple system  with a triple linear product $[\c,\c,\c]'$
defined by
$$[a\o t^i,b\o t^j, c\o t^k]'=[a,b,c]\o t^{i+j+k},$$
for all $a\o t^i,b\o t^j, c\o t^k\in T'$, where $|a\o t^i|=|a|.$
\smallskip

%\noindent{\bf Definition 2.4.}
%An ideal of a $\delta$-Jordan Lie supertriple system  is a nonzero subspace $I$ for which $[I,T,T]\subseteq I.$
% Moreover, if $[T,I,I] = 0$, then $I$ is called an abelian ideal of $T$.
%\medskip

\noindent{\bf Definition 2.5.}
Let $T$  be a $\delta$-Jordn Lie supertriple system and $k$  a nonnegative integer.
 A homogeneous linear map $D: T\rightarrow T$  is said to be a homogeneous $k$-derivation of $T$  if it satisfies
\begin{eqnarray}
\delta^{k}[D(a),b,c]+\delta^{k}(-1)^{|D||a|}[a, D(b),c]+\delta^{k}(-1)^{|D|(|a|+|b|)}[a,b,D(c)]=D([a,b,c]),
\end{eqnarray}
for all $a,b,c\in T$, where $|D|$ denotes the degree of $D$.
\smallskip

We denote by $Der(T)=\bigoplus_{k\geq 0}Der_{k}(T)$, where $Der_{k}(T)$ is the set of all homogeneous $k$-derivations of $T$.
Obviously, $Der(T)$ is a subalgebra of $End(T)$ and has a normal Lie superalgebra structure via the bracket product
\begin{eqnarray*}
[D,D']=DD'-(-1)^{|D||D'|}D'D.
\end{eqnarray*}

\noindent{\bf Definition 2.6.}
Let $T$  be a $\delta$-Jordan Lie supertriple system and $k$  a nonnegative integer.
$D\in End_{\overline{s}}(T)$ is said to be a homogeneous generalized $k$-derivation of $T$, if there  exist three endomorphisms
$D', D'',D'''\in End_{\overline{s}}(T)$ such that
\begin{eqnarray}
\delta^{k}[D(a),b,c]+\delta^{k}(-1)^{s|a|}[a, D'(b),c]+\delta^{k}(-1)^{s(|a|+|b|)}[a,b,D''(c)]=D'''([a,b,c]),
\end{eqnarray}
for all $a,b,c\in T$.
\smallskip

\noindent{\bf Definition 2.7.}
Let $T$  be a $\delta$-Jordan Lie supertriple system and $k$  a nonnegative integer.
$D\in End_{\overline{s}}(T)$ is said to be a homogeneous  $k$-quasiderivation of $T$, if there  exist an endomorphism
$D'\in End_{\overline{s}}(T)$ such that
\begin{eqnarray}
\delta^{k}[D(a),b,c]+\delta^{k}(-1)^{s|a|}[a, D(b),c]+\delta^{k}(-1)^{s(|a|+|b|)}[a,b,D(c)]=D'([a,b,c]),
\end{eqnarray}
for all $a,b,c\in T$.
\smallskip

Let $GDer(T)$ and $QDer(T)$ be the sets of homogeneous generalized $k$-derivations and of homogeneous  $k$-quasiderivations, respectively. That is,
\begin{eqnarray*}
GDer(T)=\bigoplus_{k\geq 0}GDer_{k}(T),~QDer(T)=\bigoplus_{k\geq 0}QDer_{k}(T).
\end{eqnarray*}
\smallskip

\noindent{\bf Definition 2.8.}
Let $T$  be a $\delta$-Jordan Lie supertriple system and $k$  a nonnegative integer.
The $k$-centroid of  $T$ is the space of linear transformations on $T$ given by
\begin{eqnarray}
&&C^{k}(T)=\{D\in End(T)|\delta^{k}[D(a),b,c]=\delta^{k}(-1)^{|D||a|}[a,D(b),c]\nonumber\\
&&~~~~~~~~~~~~~~~~~~~~~~~~~~~~~~~~~~~=\delta^{k}(-1)^{|D|(|a|+|b|)}[a,b,D(c)]=D([a,b,c])\}.~~
\end{eqnarray}
We denote $C(T)=\bigoplus_{k\geq 0}C^{k}(T)$  and call it the centroid of  $T$.
\medskip

\noindent{\bf Definition 2.9.}
Let $T$  be a $\delta$-Jordan Lie supertriple system.
The  quasicentroid of $T$ is the space of linear transformations on $T$ given by
\begin{eqnarray}
QC(T)=\{D\in End(T)|~D([a,b,c])=[D(a),b,c],\forall a,b,c\in T\}.
\end{eqnarray}

\noindent{\bf Remark 2.10.}
Let $T$  be a $\delta$-Jordan Lie supertriple system.
Then $QC(T)\subseteq C(T).$

For any $D\in QC(T)$ and $a,b,c\in T$, we have
\begin{eqnarray}
D([a,b,c])=[D(a),b,c]=(-1)^{|D||a|}[a,D(b),c]=(-1)^{|D|(|a|+|b|)}[a,b,D(c)].
\end{eqnarray}
In fact, by the definition of the $\delta$-Jordan Lie supertriple system, we have
\begin{eqnarray*}
D([a,b,c])&=&-\delta(-1)^{|a||b|}D([b,a,c])
=-\delta(-1)^{|a||b|}[D(b),a,c]\\
&=&\delta^{2}(-1)^{|a||b|}(-1)^{|a|(|b|+|D|)}[a,D(b),c]
=(-1)^{|D||a|}[a,D(b),c].
\end{eqnarray*}
Similarly, we have
\begin{eqnarray*}
(-1)^{|a|(|D|+|c|)}[a,b,D(c)]
&=&-(-1)^{|a||b|}[b,D(c),a]-(-1)^{|b|(|D|+|c|)}[D(c),a, b]\\
&=&-(-1)^{|b|(|a|+|D|)}D([b,c,a])-(-1)^{|b|(|D|+|c|)}D([c,a,b])\\
&=&-(-1)^{|b||D|}D((-1)^{|b||a|}[b,c,a]+(-1)^{|c||b|}[c,a,b])\\
&=&(-1)^{|b||D|}D((-1)^{|a||c|}[a,b,c]).
\end{eqnarray*}
It follows that
$[a,b,D(c)]=(-1)^{|D|(|a|+|b|)}[a,b,D(c)].$
\medskip

\noindent{\bf Definition 2.11.}
Let $T$  be a $\delta$-Jordan Lie supertriple system.
$D\in End(T)$ is said to be a central derivation  of $T$ if
\begin{eqnarray}
D([a,b,c])=[D(a),b,c]=0,
\end{eqnarray}
for all $a,b,c\in T$.
Denote the set of all central  derivations by $ZDer(T)$.
\medskip

\noindent{\bf Remark 2.12.}
Let $T$  be a $\delta$-Jordan Lie supertriple system.
Then
\begin{eqnarray*}
ZDer(T)\subseteq Der(T)\subseteq QDer(T)\subseteq GDer(T) \subseteq End(T).
\end{eqnarray*}

\noindent{\bf Definition 2.13.}
Let $T$  be a $\delta$-Jordan Lie supertriple system.
If $Z(T)=\{a\in T|~[a,b,c]=0,~\forall~b,c\in T\}$,
then  $Z(T)$ is called the center of $T$.
\medskip

\noindent{\bf Proposition 2.14.}
Let $T$  be a $\delta$-Jordan Lie supertriple system,
 then the following statements hold:

(1)~$GDer(T), QDer(T)$ and $C(T)$ are subalgebras of $End(T)$.

(2)~$ZDer(T)$ is an ideal of $Der(T)$.
\medskip

\noindent{\bf Proof.}
(1) We only prove that $GDer(T)$ is a subalgebra of $End(T)$, and similarly for cases of $ QDer(T)$ and $C(T)$.
For any $D_1\in GDer_{k}(T),D_2\in GDer_{l}(T)$ and $a,b,c\in T$, we have
\begin{eqnarray*}
&&[D_1D_2(a),b,c]\\
&=&\delta^{k} D'''_1[D_2(a),b,c]-(-1)^{|D_1|(|D_2|+|a|)}[D_2(a),D'_1(b),c]\\
  &&-(-1)^{|D_1|(|D_2|+|a|+|b|)}[D_2(a),b,D''_1(c)]\\
&=&\delta^{k} D'''_1\{\delta^{l} D'''_2([a,b,c])-(-1)^{|D_2||a|}[a, D'_2(b),c]-(-1)^{|D_2|(|a|+|b|)}[a,b,D''_2(c)]\}\\
  &&-(-1)^{|D_1|(|D_2|+|a|)}[D_2(a),D'_1(b),c]-(-1)^{|D_1|(|D_2|+|a|+|b|)}[D_2(a),b,D''_1(c)]\\
&=&\delta^{k+l}D'''_1D'''_2([a,b,c])-\delta^{l}(-1)^{|D_2||a|}D'''_1[a, D'_2(b),c]-\delta^{l}(-1)^{|D_2|(|a|+|b|)}D'''_1[a,b,D''_2(c)]\\
  &&-(-1)^{|D_1|(|D_2|+|a|)}[D_2(a),D'_1(b),c]-(-1)^{|D_1|(|D_2|+|a|+|b|)}[D_2(a),b,D''_1(c)]\\
&=&\delta^{k+l}D'''_1D'''_2([a,b,c])
   -\delta^{k+l}(-1)^{|D_2||a|}[D_1(a), D'_2(b),c]-\delta^{k+l}(-1)^{(|D_1|+|D_2|)|a|}[a, D'_1 D'_2(b),c]\\
   &&-\delta^{k+l}(-1)^{|D_2||a|+|D_1|(|a|+|b|+|D_2|)}[a, D'_2(b),D''_1(c)]
   -\delta^{k+l}(-1)^{|D_2|(|a|+|b|)}[D_1(a),b,D''_2(c)]\\
   &&-\delta^{k+l}(-1)^{|D_2|(|a|+|b|)+|D_1||a|}[a,D'_1(b),D''_2(c)]-\delta^{k+l}(-1)^{(|D_1|+|D_2|)(|a|+|b|)}[a,b,D''_1 D''_2(c)]\\
   &&-(-1)^{|D_1|(|D_2|+|a|)}[D_2(a),D'_1(b),c]-(-1)^{|D_1|(|D_2|+|a|+|b|)}[D_2(a),b,D''_1(c)].
\end{eqnarray*}
Similarly, we have
\begin{eqnarray*}
&&[D_2D_1(a),b,c]\\
&=&\delta^{k+l}D'''_2D'''_1([a,b,c])
   -\delta^{k+l}(-1)^{|D_1||a|}[D_2(a), D'_1(b),c]-\delta^{k+l}(-1)^{(|D_1|+|D_2|)|a|}[a, D'_2 D'_1(b),c]\\
   &&-\delta^{k+l}(-1)^{|D_1||a|+|D_2|(|a|+|b|+|D_1|)}[a, D'_1(b),D''_2(c)]
   -\delta^{k+l}(-1)^{|D_1|(|a|+|b|)}[D_2(a),b,D''_1(c)]\\
   &&-\delta^{k+l}(-1)^{|D_1|(|a|+|b|)+|D_2||a|}[a,D'_2(b),D''_1(c)]-\delta^{k+l}(-1)^{(|D_1|+|D_2|)(|a|+|b|)}[a,b,D''_2 D''_1(c)]\\
   &&-(-1)^{|D_2|(|D_1|+|a|)}[D_1(a),D'_2(b),c]-(-1)^{|D_2|(|D_1|+|a|+|b|)}[D_1(a),b,D''_2(c)].
\end{eqnarray*}
It follows that
\begin{eqnarray*}
[[D_1, D_2](a),b,c]
&=&[D_1D_2(a),b,c]-(-1)^{k+l}[D_2D_1(a),b,c]\\
&=&\delta^{k+l}(D'''_1D'''_2-(-1)^{|D_1|+|D_2|}D'''_2D'''_1)[a,b,c]\\
&&-(-1)^{|a|(|D_1|+|D_2|)}[a,(D'_1D'_2-(-1)^{|D_1|+|D_2|}D'_2D'_1)(b),c]\\
&&-(-1)^{(|a|+|b|)(|D_1|+|D_2|)}[a,b,(D''_1D''_2-(-1)^{|D_1|+|D_2|}D''_2D''_1)(c)]
\end{eqnarray*}
\begin{eqnarray*}
&=&\delta^{k+l}[D'''_1,D'''_2][a,b,c]
   -(-1)^{|a|(|D_1|+|D_2|)}[a,[D'_1,D'_2](b),c]\\
   &&-(-1)^{(|a|+|b|)(|D_1|+|D_2|)}[a,b,[D''_1,D''_2](c)].
\end{eqnarray*}
Obviously, $[D'_1,D'_2],[D''_1,D''_2]$ and $[D'''_1,D'''_2]$ are contained in $End(T)$,
thus $[D_1, D_2]\in GDer_{k+l}(T)\subseteq GDer(T)$, that is, $GDer(T)$ is a subalgebra of $End(T)$.
\medskip

(2) For any $D_1\in ZDer(T),D_2\in Der_{k}(T)$ and $a,b,c\in T$, we have
\begin{eqnarray*}
[D_1,D_2]([a,b,c])
=D_1D_2([a,b,c])-(-1)^{|D_1||D_2|}D_2D_1([a,b,c])=0.
\end{eqnarray*}
Also, we have
\begin{eqnarray*}
[[D_1,D_2](a),b,c]
&=&[D_1D_2(a),b,c]-(-1)^{|D_1||D_2|}[D_2D_1(a),b,c]\\
&=&0-(-1)^{|D_1||D_2|}[D_2D_1(a),b,c]\\
&=&-(-1)^{|D_1||D_2|}(\delta^{k}D_2([D_1(a),b,c])-(-1)^{|D_2|(|D_1|+|a|)}[D_1(a),D_2(b),c])\\
&&+(-1)^{|D_1||D_2|}(-1)^{|D_2|(|D_1|+|a|+|b|)}[D_1(a),b,D_2(c)]\\
&=&0.
\end{eqnarray*}
It follows that $[D_1,D_2]\in ZDer(T)$.
That is, $ZDer(T)$ is an ideal of $Der(T)$. $\hfill \Box$
\medskip

\noindent{\bf Proposition 2.15.}
Let $T$  be a $\delta$-Jordan Lie supertriple system,
 then the following statements hold:

(1)~$[Der(T), C(T)]\subseteq C(T)$.

(2)~$[QDer(T), QC(T)]\subseteq QC(T)$.

(3)~$[QC(T), QC(T)]\subseteq QDer(T)$.

(4)~$C(T)\subseteq QDer(T)$.
\medskip

\noindent{\bf Proof.}
(1)~For any $D_1\in Der_{k}(T), D_2\in C_{l}(T)$ and $a,b,c\in T$, we have
\begin{eqnarray*}
&&[D_1,D_2]([a,b,c])\\
&=&D_1D_2([a,b,c])-(-1)^{|D_1||D_2|}D_2D_1([a,b,c])\\
&=&\delta^{l}D_1([D_2(a),b,c])-\delta^{k}(-1)^{|D_1||D_2|}D_2([D_1(a),b,c])\\
   &&-\delta^{k}(-1)^{|D_1|(|D_2|+|a|)}D_2([a,D_1(b),c])-\delta^{k}(-1)^{|D_1|(|D_2|+|a|+|b|)}D_2([a,b,D_1(c)])\\
&=&\delta^{k+l}[D_1 D_2(a),b,c]+\delta^{k+l}(-1)^{|D_1|(|D_2|+|a|)}[ D_2(a),D_1(b),c]\\
   &&+\delta^{k+l}(-1)^{|D_1|(|D_2|+|a|+|b|)}[ D_2(a),b,D_1(c)]-\delta^{k+l}(-1)^{|D_1||D_2|}[D_2 D_1(a),b,c]\\
   &&-\delta^{k+l}(-1)^{|D_1|(|D_2|+|a|)}[D_2(a),D_1(b),c]-\delta^{k+l}(-1)^{|D_1|(|D_2|+|a|+|b|)}[D_2(a),b,D_1(c)]\\
&=&\delta^{k+l}[[D_1,D_2](a),b,c].
\end{eqnarray*}
Similarly, one can check that
$$[[D_1,D_2](a),b,c]=(-1)^{|a|(k+l)}[a,[D_1,D_2](b),c]
=(-1)^{(|a|+|b|)(k+l)}[a,b,[D_1,D_2](c)].$$
It follows that $[D_1,D_2]\in C(T)$, thus $[Der(T), C(T)]\subseteq C(T)$.
\medskip

(2) Similar to the proof of (1).
\medskip

(3) ~For any $D_1, D_2\in QC(T)$ and $a,b,c\in T$, we have
\begin{eqnarray*}
&&[[D_1,D_2](a),b,c]+(-1)^{|a|(|D_1|+|D_2|)}[a,[D_1,D_2](b),c]\\
&&~~+(-1)^{(|a|+|b|)(|D_1|+|D_2|)}[a,b,[D_1,D_2](c)]\\
&=&[D_1D_2(a),b,c]+(-1)^{|a|(|D_1|+|D_2|)}[a,D_1D_2(b),c]\\
&&~~+(-1)^{(|a|+|b|)(|D_1|+|D_2|)}[a,b,D_1D_2(c)]-(-1)^{|D_1||D_2|}[D_2D_1(a),b,c]\\
&&~~-(-1)^{|a|(|D_1|+|D_2|)+|D_1||D_2|}[a,D_2D_1(b),c]-(-1)^{(|a|+|b|)(|D_1|+|D_2|)+|D_1||D_2|}[a,b,D_2D_1(c)].
\end{eqnarray*}
Since $D_1, D_2\in QC(T)$, we have
\begin{eqnarray*}
[D_1D_2(a),b,c]
&=&(-1)^{|D_1|(|a|+|D_2|)}[D_2(a),D_1(b),c]\\
&=&(-1)^{|D_1|(|a|+|D_2|)+|a||D_2|}[a,D_2D_1(b),c].
\end{eqnarray*}
Similarly, one may check that
\begin{eqnarray*}
&&[a,D_1D_2(b),c]=(-1)^{|b|(|D_1|+|D_2|)+|D_1||D_2|}[a,b,D_2D_1(c)],\\
&&[a,b,D_1D_2(c)]=(-1)^{(|a|+|b|)(|D_1|+|D_2|)+|D_1||D_2|}[D_2D_1(a),b,c].
\end{eqnarray*}
It follows that
\begin{eqnarray*}
&&[[D_1,D_2](a),b,c]+(-1)^{|a|(|D_1|+|D_2|)}[a,[D_1,D_2](b),c]\\
&&~~~~~~~~~~~~~~~~~~~~~~+(-1)^{(|a|+|b|)(|D_1|+|D_2|)}[a,b,[D_1,D_2](c)]=0.
\end{eqnarray*}
Therefore, $[D_1,D_2]\subseteq QDer(T)$  and $[QC(T), QC(T)]\subseteq QDer(T)$.
\medskip

(4) ~For any $D\in C_{k}(T)$ and $a,b,c\in T$, we have
\begin{eqnarray*}
D([a,b,c])
=\delta^{k}[D(a),b,c]
=\delta^{k}(-1)^{|D||a|}[a, D(b),c]
=\delta^{k}(-1)^{|D|(|a|+|b|)}[a,b,D(c)].
\end{eqnarray*}
Thus
$$
\delta^{k}[D(a),b,c]
+\delta^{k}(-1)^{|D||a|}[a, D(b),c]
+\delta^{k}(-1)^{|D|(|a|+|b|)}[a,b,D(c)]=
3D([a,b,c]),
$$
that is, $D\in QDer_{k}(T)$  and $QC(T)\subseteq QDer(T)$.
$\hfill \Box$
\medskip

\noindent{\bf Theorem 2.16.}
Let $T$  be a $\delta$-Jordan Lie supertriple system,
 then $[C(T), QC(T)]\subseteq End(T,Z(T))$.
  Moreover, if $Z(T)=\{0\}$, then  $[C(T), QC(T)]=\{0\}$.
\medskip

\noindent{\bf Proof.}
For any $D_1\in C_{k}(T), D_2\in QC_{l}(T)$ and $a,b,c\in T$, we have
\begin{eqnarray*}
[[D_1,D_2](a),b,c]
&=&[D_1D_2(a),b,c]-(-1)^{|D_1||D_2|}[D_2D_1(a),b,c]\\
&=&\delta^{k}D_1([D_2(a),b,c])-(-1)^{|D_1||D_2|}(-1)^{(|D_1|+|a|)|D_2|}[D_1(a),D_2(b),c]\\
&=&\delta^{k}(-1)^{|a||D_2|}D_1([a,D_2(b),c])-(-1)^{|a||D_2|}D_1([a,D_2(b),c])\\
&=&0.
\end{eqnarray*}
So $[D_1,D_2](a)\subseteq Z(T)$ and therefore $[C(T), QC(T)]\subseteq End(T,Z(T))$.
 Moreover, if $Z(T)=\{0\}$, then  it is easy to see that $[C(T), QC(T)]=\{0\}$.
$\hfill \Box$
\medskip

\noindent{\bf Theorem 2.17.}
Let $T$  be a $\delta$-Jordan Lie supertriple system and $Char~k\neq 2$,
 then $ZDert=C(T)\cap Der(T)$.
\medskip

\noindent{\bf Proof.}
For any $D\in C_{k}(T)\cap Der_{k}(T)$ and $a,b,c\in T$, we have
\begin{eqnarray*}
&&D([a,b,c])=\delta^{k}[D(a),b,c]+\delta^{k}(-1)^{|D||a|}[a, D(b),c]+\delta^{k}(-1)^{|D|(|a|+|b|)}[a,b,D(c)],\\
&&D([a,b,c])=\delta^{k}[D(a),b,c]=\delta^{k}(-1)^{|D||a|}[a, D(b),c]=\delta^{k}(-1)^{|D|(|a|+|b|)}[a,b,D(c)].
\end{eqnarray*}
It follows that $2D([a,b,c])=0$.
Thus $D([a,b,c])=0$ since $Char~k\neq 2$,
that is, $D\in ZDert$. So $C(T)\cap Der(T)\in ZDert$.

On the other hand, for any $D\in ZDert$ and $a,b,c\in T$, we have $D([a,b,c])=0$.
Clearly, Eq. (2.5) and Eq. (2.8) hold, that is,  $D\in C_{k}(T)\cap Der_{k}(T)$ and therefore  $ ZDert\in C(T)\cap Der(T)$.
And this completes the proof.
$\hfill \Box$

\section{The  cohomology of $\delta$-Jordan Lie supertriple systems}
\def\theequation{\arabic{section}. \arabic{equation}}
\setcounter{equation} {0}

In this section, we introduce the notion of the representation of $\delta$-Jordan Lie supertriple systems
and present its low-dimensional cohomologies.
\medskip

\noindent{\bf Definition 3.1.}
Let $T$  be a $\delta$-Jordan Lie supertriple system and $V$  a $Z_2$-graded  vector space.
Suppose that there exists a bilinear mapping $\theta:T\o T\rightarrow End(V)$  satisfying the following axioms:
\begin{eqnarray}
&&(-1)^{(|a|+|b|)(|c|+|d|)}\theta(c,d)\theta(a,b)-\delta(-1)^{(|a||b|+|d|(|c|+|a|))}\theta(b,d)\theta(a,c)\nonumber\\
   &&~~~~~~~~-\theta(a,[b,c,d])+(-1)^{|a|(|b|+|c|)}D(b,c)\theta(a,d)=0,\\
&&\delta(-1)^{(|a|+|b|)(|c|+|d|)}\theta(c,d)D(a,b)-\delta D(a,b)\theta(c,d)\nonumber\\
   &&~~~~~~~~+\theta([a,b,c],d)+\delta(-1)^{|c|(|a|+|b|)}\theta(c,[a,b,d])=0,\\
&&D([a,b,c],d)+(-1)^{|c|(|a|+|b|)}D(c,[a,b,d])\nonumber\\
 &&~~~~~~~~-\delta D(a,b)D(c,d)+(-1)^{(|a|+|b|)(|c|+|d|)}D(c,d) D(a,b),
\end{eqnarray}
for $a,b,c,d\in T$, where $D(a,b)=(-1)^{|a||b|}\theta(b,a)-\delta\theta(a,b)$,
 then $(V,\theta)$ is called the representation of $T$, $V$ is called a $T$-module.
\medskip

\noindent{\bf Example 3.2.}
Let $T$  be a $\delta$-Jordan Lie supertriple system.
Define  $\theta:T\o T\rightarrow End(T)$ by $\theta(a,b)(x)=(-1)^{|x|(|a|+|b|)}[x,a,b]$.
It is not hard to check that $D(a,b)(x)=\delta[a,b,x]$ and $T$ itself is a $T$-module.
In this case,  $(T,\theta)$ is said to be the adjoint representation of $T$.
\medskip

\noindent{\bf Proposition 3.3.}
Let $T$  be a $\delta$-Jordan Lie supertriple system and $(V,\theta)$ the representation.
Then $T\oplus V$ has a structure of a $\delta$-Jordan Lie supertriple system.
\medskip

\noindent{\bf Proof.}
Define a  triple linear product $[\c, \c, \c]:(T\oplus V)\o(T\oplus V)\o(T\oplus V)\rightarrow (T\oplus V)$ by
\begin{eqnarray*}
&&[(a,u),(b,v),(c,w)]\\
&&~~~~~=([a,b,c],(-1)^{|a|(|b|+|c|)}\theta(b,c)(u)-\delta(-1)^{|b||c|}\theta(a,c)(v)+\delta D(a,b)(w)),
\end{eqnarray*}
for all $(a,u),(b,v),(c,w)\in T\oplus V$, where $|(a,u)|=|a|$.
\smallskip

Now we check that the operation $[\c, \c, \c]$ defined above satisfies axioms in Definition 2.1.
It is easy to see that  Eq. (2.1) holds since $T$  is a $\delta$-Jordan Lie supertriple system.
\smallskip

For Eq. (2.2), we take any $(a,u),(b,v),(c,w)\in T\oplus V$ and compute
 \begin{eqnarray*}
&&-\delta(-1)^{|a||b|}[(b,v),(a,u),(c,w)]\\
&&~~~~~=-\delta(-1)^{|a||b|}([b,a,c],(-1)^{|b|(|a|+|c|)}\theta(a,c)(v)-\delta(-1)^{|a||c|}\theta(b,c)(u)+\delta D(b,a)(w))\\
&&~~~~~=(-\delta(-1)^{|a||b|}[b,a,c],-\delta(-1)^{|b||c|}\theta(a,c)(v)-(-1)^{|a|(|b|+|c|)}\theta(b,c)(u)\\
      &&~~~~~~~~~-\delta(-1)^{|a||b|}\delta D(b,a)(w))\\
&&~~~~~= ([a,b,c],(-1)^{|a|(|b|+|c|)}\theta(b,c)(u)-\delta(-1)^{|b||c|}\theta(a,c)(v)+\delta D(a,b)(w))\\
&&~~~~~=[(a,u),(b,v),(c,w)].
\end{eqnarray*}
The last equality holds since
$
-(-1)^{|a||b|}[b,a,w]
=[a,b,w]
=\delta  D(a,b)(w).$
\smallskip

For Eq. (2.3), we have
 \begin{eqnarray*}
&&(-1)^{|a||c|}[(a,u),(b,v),(c,w)]+(-1)^{|b||a|}[(b,v),(c,w),(a,u)]\\
&&~~~~+(-1)^{|c||b|}[(c,w),(a,u),(b,v)]\\
&=&((-1)^{|a||c|}[a,b,c]+(-1)^{|b||a|}[b,c,a]+(-1)^{|c||b|}[c,a,b],\Omega)\\
&=&(0,\Omega),
\end{eqnarray*}
where
 \begin{eqnarray*}
\Omega
&=&(-1)^{|a||b|}\theta(b,c)(u)-\delta(-1)^{(|a|+|b|)|c|}\theta(a,c)(v)+\delta(-1)^{|a||c|}D(a,b)(w)\\
&&+(-1)^{|b||c|}\theta(c,a)(v)-\delta(-1)^{(|b|+|c|)|a|}\theta(b,a)(w)+\delta(-1)^{|b||a|}D(b,c)(u)\\
&&+(-1)^{|c||b|}\theta(a,b)(w)-\delta(-1)^{(|c|+|a|)|b|}\theta(c,b)(u)+\delta(-1)^{|c||b|}D(c,a)(v)\\
&=&(-1)^{|a||b|}\theta(b,c)(u)+(-1)^{|b||c|}\theta(c,a)(v)+(-1)^{|c||b|}\theta(a,b)(w)\\
&&-(-1)^{|a||c|}\theta(a,b)(w)+(-1)^{|b||a|}\theta(b,c)(u)-(-1)^{|b||c|}\theta(c,a)(v)\\
&=&0.
\end{eqnarray*}
The second equality holds since $D(a,b)=(-1)^{|a||b|}\theta(b,a)-\delta\theta(a,b)$.
Then we have
 \begin{eqnarray*}
&&(-1)^{|a||c|}[(a,u),(b,v),(c,w)]+(-1)^{|b||a|}[(b,v),(c,w),(a,u)]\\
&&~~~~+(-1)^{|c||b|}[(c,w),(a,u),(b,v)]=(0,0),
\end{eqnarray*}
as desired.
\smallskip

For Eq. (2.4), we take any $(a,u),(b,v),(c,w),(d,m),(e,n)\in T\oplus V$.
First, we calculate the following expression:
\begin{eqnarray*}
&&[[(a,u),(b,v),(c,w)],(d,m),(e,n)]\\
&=&[([a,b,c],(-1)^{|a|(|b|+|c|)}\theta(b,c)(u)-\delta(-1)^{|b||c|}\theta(a,c)(v)+\delta D(a,b)(w)),(d,m),(e,n)]\\
&=&([[a,b,c],d,e],\Omega_1),
\end{eqnarray*}
where
\begin{eqnarray*}
\Omega_1
&=&(-1)^{(|a|+|b|+|c|)(|d|+|e|)}(-1)^{|a|(|b|+|c|)}\theta(d,e)\theta(b,c)(u)\\
 &&-\delta(-1)^{(|a|+|b|+|c|)(|d|+|e|)}(-1)^{|b||c|}\theta(d,e)\theta(a,c)(v)\\
 &&+\delta(-1)^{(|a|+|b|+|c|)(|d|+|e|)}\theta(d,e)D(a,b)(w)\\
 &&-\delta(-1)^{|d||e|}\theta([a,b,c],e)(m)
   +\delta D([a,b,c],d)(n).
\end{eqnarray*}
Second, we compute the expression $(-1)^{|c|(|a|+|b|)}[(c,w),[(a,u),(b,v),(d,m)],(e,n)]$:
\begin{eqnarray*}
&&(-1)^{|c|(|a|+|b|)}[(c,w),[(a,u),(b,v),(d,m)],(e,n)]\\
&=&(-1)^{|c|(|a|+|b|)}[(c,w),([a,b,d],(-1)^{|a|(|b|+|d|)}\theta(b,d)(u)-\delta(-1)^{|b||d|}\theta(a,d)(v)\\
     &&+\delta D(a,b)(m)),(e,n)]\\
&=&((-1)^{|c|(|a|+|b|)}[c,[a,b,d],e],\Omega_2),
\end{eqnarray*}
where
\begin{eqnarray*}
\Omega_2
&=&-\delta(-1)^{(|a|+|b|+|d|)|e|}(-1)^{|c|(|a|+|b|)+|a|(|b|+|d|)}\theta(c,e)\theta(b,d)(u)\\
 &&+(-1)^{(|a|+|b|+|d|)|e|}(-1)^{|c|(|a|+|b|)+|b||d|}\theta(c,e)\theta(a,d)(v)\\
 &&+(-1)^{(|a|+|b|+|d|+|e|)|c|}(-1)^{|c|(|a|+|b|)}\theta([a,b,d],e)(w)\\
 &&-(-1)^{(|a|+|b|+|d|)|e|}(-1)^{|c|(|a|+|b|)}\theta(c,e)D(a,b)(m)\\
  && +\delta(-1)^{|c|(|a|+|b|)} D(c,[a,b,d])(n).
\end{eqnarray*}
Third, we compute the expression $\delta(-1)^{(|a|+|b|)(|c|+|d|)}[(c,w),(d,m),[(a,u),(b,v),(e,n)]]$:
\begin{eqnarray*}
&&\delta(-1)^{(|a|+|b|)(|c|+|d|)}[(c,w),(d,m),[(a,u),(b,v),(e,n)]]\\
&=&\delta(-1)^{(|a|+|b|)(|c|+|d|)}[(c,w),(d,m),([a,b,e],(-1)^{|a|(|b|+|e|)}\theta(b,e)(u)\\
     &&-\delta(-1)^{|b||e|}\theta(a,e)(v)+\delta D(a,b)(n))]\\
&=&(\delta(-1)^{(|a|+|b|)(|c|+|d|)}[c,d,[a,b,e]],\Omega_3),
\end{eqnarray*}
where
\begin{eqnarray*}
\Omega_3
&=&(-1)^{(|a|+|b|)(|c|+|d|)}(-1)^{|a|(|b|+|e|)}D(c,d)\theta(b,e)(u)\\
 &&-\delta(-1)^{(|a|+|b|)(|c|+|d|)}(-1)^{|b||e|}D(c,d)\theta(a,e)(v)\\
 &&+\delta(-1)^{(|a|+|b|)(|c|+|d|)}(-1)^{(|a|+|b|+|d|+|e|)|c|}\theta(d,[a,b,e])(w)\\
 &&-(-1)^{(|a|+|b|)(|c|+|d|)}(-1)^{(|a|+|b|+|e|)|d|}\theta(c,[a,b,e])(m)\\
  && +\delta(-1)^{(|a|+|b|)(|c|+|d|)} D(c,d) D(a,b)(n).
\end{eqnarray*}
Four, we compute the expression $[(a,u),(b,v),[(c,w),(d,m),(e,n)]]$:
\begin{eqnarray*}
&&[(a,u),(b,v),[(c,w),(d,m),(e,n)]]\\
&=&[(a,u),(b,v),([c,d,e],(-1)^{|c|(|d|+|e|)}\theta(d,e)(w)-\delta(-1)^{|d||e|}\theta(c,e)(m)+\delta D(c,d)(n))]\\
&=&([a,b,[c,d,e]],\Omega_4),
\end{eqnarray*}
where
\begin{eqnarray*}
\Omega_4
&=&(-1)^{|a|(|b|+|c|+|d|+|e|)}\theta(b,[c,d,e])(u)
 -\delta(-1)^{|b|(|c|+|d|+|e|)}\theta(a,[c,d,e])(v)\\
 &&+\delta(-1)^{|c|(|d|+|e|)}D(a,b)\theta(d,e)(w)
 -(-1)^{|d||e|}D(a,b)\theta(c,e)(m)
  + D(a,b)D(c,d)(n).
\end{eqnarray*}
Finally, by Eq. (3.1), Eq. (3.2) and Eq. (3.3), we have
\begin{eqnarray*}
&&[[(a,u),(b,v),(c,w)],(d,m),(e,n)]
+(-1)^{|c|(|a|+|b|)}[(c,w),[(a,u),(b,v),(d,m)],(e,n)]\\
&&+\delta(-1)^{(|a|+|b|)(|c|+|d|)}[(c,w),(d,m),[(a,u),(b,v),(e,n)]]\\
&=&([[a,b,c],d,e]+(-1)^{|c|(|a|+|b|)}[c,[a,b,d],e]+\delta(-1)^{(|a|+|b|)(|c|+|d|)}[c,d,[a,b,e]],\Omega_1+\Omega_2+\Omega_3)\\
&=&([a,b,[c,d,e]],\Omega_4)\\
&=&[(a,u),(b,v),[(c,w),(d,m),(e,n)]],
\end{eqnarray*}
as desired, and this finishes the proof.
$\hfill \Box$
\medskip

\noindent{\bf Corollary 3.4.}
 Any $\delta$-Jordan Lie supertriple system   can be considered as a subspace of a $\delta$-Jordan Lie superalgebra.
\medskip

\noindent{\bf Proof.}
Straightforward from Example 3.2 and Proposition 3.3.
$\hfill \Box$
\medskip

\noindent{\bf Definition 3.5.}
Let $T$  be a $\delta$-Jordan Lie supertriple system and $V$  a $T$-module by a bilinear map $\theta$.
If an $n$-linear map $f:T\times T\times\cdots \times T\rightarrow T$ satisfies the following axioms:

(1)~$f(x_1,x_2,\cdots,x,y,\cdots,x_n)=-\delta(-1)^{|x||y|}f(x_1,x_2,\cdots,x,y,\cdots,x_n)$;

(2)~$(-1)^{|x||z|}f(x_1,x_2,\cdots,x_{n-3},x,y,z)+(-1)^{|y||x|}f(x_1,x_2,\cdots,x_{n-3},y,z,x)\\
~~~~~~~~~~~~~~~+(-1)^{|z||y|}f(x_1,x_2,\cdots,x_{n-3},z,x,y)=0,$\\
then $f$ is called an $n$-cochain on $T$. Denote by $C_{\delta}^{n}(T,V)$  the set of all $n$-cochains, for $n\geq 1$.
\medskip

\noindent{\bf Definition 3.6.}
Let $T$  be a $\delta$-Jordan Lie supertriple system and $V$  a $T$-module by a bilinear map $\theta$.
For $n=1,2$, the coboundary operator $d_n: C_{\delta}^{n}(T,V)\rightarrow C_{\delta}^{n+2}(T,V)$  is  defined as follow:

$\bullet$ If $f\in C_{\delta}^{1}(T,V)$, then
\begin{eqnarray*}
d^{1}f(x_1,x_2,x_3)
&=&(-1)^{(|f|+|x_1|)(|x_2|+|x_3|)}\theta(x_2,x_3)f(x_1)-f([x_1,x_2,x_3])\\
    &&+\delta(-1)^{|f|(|x_1|+|x_2|)}D(x_1,x_2)f(x_3)\\
    && -\delta(-1)^{|x_2||x_3|+|f|(|x_1|+|x_3|)}\theta(x_1,x_3)f(x_2).
\end{eqnarray*}

$\bullet$ If $f\in C_{\delta}^{2}(T,V)$, then
\begin{eqnarray*}
d^{2}f(y,x_1,x_2,x_3)
&=&(-1)^{(|f|+|y|+|x_1|)(|x_2|+|x_3|)}\theta(x_2,x_3)f(y,x_1) -f(y,[x_1,x_2,x_3])\\
  && -\delta(-1)^{|x_2||x_3|+(|f|+|y_1|)(|x_1|+|x_3|)}\theta(x_1,x_3)f(y,x_2)\\
  &&+\delta(-1)^{(|f|+|y|)(|x_1|+|x_2|)}D(x_1,x_2)f(y,x_3).
\end{eqnarray*}

$\bullet$ If $f\in C_{\delta}^{3}(T,V)$, then
\begin{eqnarray*}
&&d^{3}f(x_1,x_2,x_3,x_4,x_5)\\
&&~~~~=(-1)^{(|f|+|x_1|+|x_2|+|x_3|)(|x_4|+|x_5|)}\theta(x_4,x_5)f(x_1,x_2,x_3)\\
  &&~~~~~~~~-\delta(-1)^{(|f|+|x_1|+|x_2|)(|x_3|+|x_5|)+|x_4||x_5|}\theta(x_3,x_5)f(x_1,x_2,x_4)\\
  &&~~~~~~~~-\delta(-1)^{|f|(|x_1|+|x_2|)}D(x_1,x_2)f(x_3,x_4,x_5)\\
  &&~~~~~~~~+(-1)^{(|f|+|x_1|+|x_2|)(|x_3|+|x_4|)}D(x_3,x_4)f(x_1,x_2,x_5)\\
  &&~~~~~~~~+f([x_1,x_2,x_3],x_4,x_5)-f(x_1,x_2,[x_3,x_4,x_5])\\
  &&~~~~~~~~+(-1)^{|x_3|(|x_1|+|x_2|)}f(x_3,[x_1,x_2,x_4],x_5)\\
  &&~~~~~~~~+\delta(-1)^{(|x_1|+|x_2|)(|x_3|+|x_4|)}f(x_3,x_4,[x_1,x_2,x_5]).
\end{eqnarray*}

$\bullet$ If $f\in C_{\delta}^{4}(T,V)$, then
\begin{eqnarray*}
&&d^{4}f(y,x_1,x_2,x_3,x_4,x_5)\\
&&~~~~=(-1)^{(|f|+|y|+|x_1|+|x_2|+|x_3|)(|x_4|+|x_5|)}\theta(x_4,x_5)f(y,x_1,x_2,x_3)\\
&&~~~~~~~~-\delta(-1)^{(|f|+|y|+|x_1|+|x_2|)(|x_3|+|x_5|)+|x_4||x_5|}\theta(x_3,x_5)f(y,x_1,x_2,x_4)\\
&&~~~~~~~~-\delta(-1)^{(|f|+|y|)(|x_1|+|x_2|)}D(x_1,x_2)f(y,x_3,x_4,x_5)\\
 &&~~~~~~~~+(-1)^{(|f|+|y|+|x_1|+|x_2|)(|x_3|+|x_4|)}D(x_3,x_4)f(y,x_1,x_2,x_5)\\
  &&~~~~~~~~+f(y,[x_1,x_2,x_3],x_4,x_5)-f(y,x_1,x_2,[x_3,x_4,x_5])\\
    &&~~~~~~~~+(-1)^{|x_3|(|x_1|+|x_2|)}f(y,x_3,[x_1,x_2,x_4],x_5)\\
    &&~~~~~~~~+\delta(-1)^{(|x_1|+|x_2|)(|x_3|+|x_4|)}f(y,x_3,x_4,[x_1,x_2,x_5]).
\end{eqnarray*}

\noindent{\bf Theorem 3.7.}
Let $T$  be a $\delta$-Jordan Lie supertriple system and $V$  a $T$-module by a bilinear map $\theta$.
The coboundary operator $d^n$ defined above satisfies $d^{n+2}d^{n}=0$, $n=1,2$.

\noindent{\bf Proof.}
From the definition of the coboundary operator, it follows immediately that
$d^{3}d^{1}=0$ implies $d^{4}d^{2}=0$ . Then, we only need to prove $d^{3}d^{1}=0$.
\begin{eqnarray}
&&d^{3}(d^{1}f)(x_1,x_2,x_3,x_4,x_5)\nonumber\\
&&~~~~~~=(-1)^{(|f|+|x_1|+|x_2|+|x_3|)(|x_4|+|x_5|)}\theta(x_4,x_5)(d^{1}f)(x_1,x_2,x_3)\\
  &&~~~~~~~~-\delta(-1)^{(|f|+|x_1|+|x_2|)(|x_3|+|x_5|)+|x_4||x_5|}\theta(x_3,x_5)(d^{1}f)(x_1,x_2,x_4)\\
  &&~~~~~~~~-\delta(-1)^{|f|(|x_1|+|x_2|)}D(x_1,x_2)(d^{1}f)(x_3,x_4,x_5)\\
  &&~~~~~~~~+(-1)^{(|f|+|x_1|+|x_2|)(|x_3|+|x_4|)}D(x_3,x_4)(d^{1}f)(x_1,x_2,x_5)\\
  &&~~~~~~~~+(d^{1}f)([x_1,x_2,x_3],x_4,x_5)-(d^{1}f)(x_1,x_2,[x_3,x_4,x_5])\\
  &&~~~~~~~~+(-1)^{|x_3|(|x_1|+|x_2|)}(d^{1}f)(x_3,[x_1,x_2,x_4],x_5)\\
  &&~~~~~~~~+\delta(-1)^{(|x_1|+|x_2|)(|x_3|+|x_4|)}(d^{1}f)(x_3,x_4,[x_1,x_2,x_5]).
\end{eqnarray}
By Definition 3.6, we have
\begin{eqnarray}
(3.4)&=&(-1)^{(|f|+|x_1|+|x_2|+|x_3|)(|x_4|+|x_5|)}\theta(x_4,x_5)(d^{1}f)(x_1,x_2,x_3)\nonumber\\
     &=&(-1)^{(|f|+|x_1|+|x_2|+|x_3|)(|x_4|+|x_5|)}\theta(x_4,x_5)
          \{-f([x_1,x_2,x_3])\nonumber\\
             &&+(-1)^{(|f|+|x_1|)(|x_2|+|x_3|)}\theta(x_2,x_3)f(x_1)+\delta(-1)^{|f|(|x_1|+|x_2|)}D(x_1,x_2)f(x_3)\nonumber\\
             && -\delta(-1)^{|x_2||x_3|+|f|(|x_1|+|x_3|)}\theta(x_1,x_3)f(x_2)\},\\
(3.5)&=& -\delta(-1)^{(|f|+|x_1|+|x_2|)(|x_3|+|x_5|)+|x_4||x_5|}\theta(x_3,x_5)(d^{1}f)(x_1,x_2,x_4)\nonumber\\
     &=&-\delta(-1)^{(|f|+|x_1|+|x_2|)(|x_3|+|x_5|)+|x_4||x_5|}\theta(x_3,x_5)
          \{-f([x_1,x_2,x_4])\nonumber\\
             &&+(-1)^{(|f|+|x_1|)(|x_2|+|x_4|)}\theta(x_2,x_4)f(x_1)+\delta(-1)^{|f|(|x_1|+|x_2|)}D(x_1,x_2)f(x_4)\nonumber\\
             && -\delta(-1)^{|x_2||x_4|+|f|(|x_1|+|x_4|)}\theta(x_1,x_4)f(x_2)\},\\
(3.6)&=&-\delta(-1)^{|f|(|x_1|+|x_2|)}D(x_1,x_2)(d^{1}f)(x_3,x_4,x_5) \nonumber\\
      &=&-\delta(-1)^{|f|(|x_1|+|x_2|)}D(x_1,x_2)
      \{-f([x_3,x_4,x_5])\nonumber\\
             &&+(-1)^{(|f|+|x_3|)(|x_4|+|x_5|)}\theta(x_4,x_5)f(x_3)+\delta(-1)^{|f|(|x_3|+|x_4|)}D(x_3,x_4)f(x_5)\nonumber\\
             && -\delta(-1)^{|x_4||x_5|+|f|(|x_3|+|x_5|)}\theta(x_3,x_5)f(x_4)\},\\
(3.7)&=&(-1)^{(|f|+|x_1|+|x_2|)(|x_3|+|x_4|)}D(x_3,x_4)(d^{1}f)(x_1,x_2,x_5)  \nonumber\\
       &=&(-1)^{(|f|+|x_1|+|x_2|)(|x_3|+|x_4|)}D(x_3,x_4)
        \{-f([x_1,x_2,x_5])\nonumber\\
             &&+(-1)^{(|f|+|x_1|)(|x_2|+|x_5|)}\theta(x_2,x_5)f(x_1)+\delta(-1)^{|f|(|x_1|+|x_2|)}D(x_1,x_2)f(x_5)\nonumber\\
             && -\delta(-1)^{|x_2||x_5|+|f|(|x_1|+|x_5|)}\theta(x_1,x_5)f(x_2)\},\\
(3.8)&=&(d^{1}f)([x_1,x_2,x_3],x_4,x_5)-(d^{1}f)(x_1,x_2,[x_3,x_4,x_5])\nonumber\\
     &=&(-1)^{(|f|+|x_1|+|x_2|+|x_3|)(|x_4|+|x_5|)}\theta(x_4,x_5)f([x_1,x_2,x_3])-f([[x_1,x_2,x_3],x_4,x_5])\nonumber\\
      &&+\delta(-1)^{|f|(|x_1|+|x_2|+|x_3|+|x_4|)}D([x_1,x_2,x_3],x_4)f(x_5)\nonumber\\
             && -\delta(-1)^{|x_4||x_5|+|f|(|x_1|+|x_2|+|x_3|+|x_5|)}\theta([x_1,x_2,x_3],x_5)f(x_4)\nonumber
 \end{eqnarray}
\begin{eqnarray}
             && -(-1)^{(|f|+|x_1|)(|x_2|+|x_3|+|x_4|+|x_5|)}\theta(x_2,[x_3,x_4,x_5])f(x_1)\nonumber\\
             &&+f([x_1,x_2,[x_3,x_4,x_5]])-\delta(-1)^{|f|(|x_1|+|x_2|)}D(x_1,x_2)f([x_3,x_4,x_5])\nonumber\\
             &&+\delta(-1)^{|f|(|x_1|+|x_3|+|x_4|+|x_5|)+|x_2|(|x_3|+|x_4|+|x_5|)}\theta(x_1,[x_3,x_4,x_5],)f(x_2),\\
(3.9)&=&(-1)^{|x_3|(|x_1|+|x_2|)}(d^{1}f)(x_3,[x_1,x_2,x_4],x_5) \nonumber\\
      &=&(-1)^{|x_3|(|x_1|+|x_2|)+(|f|+|x_3|)(|x_1|+|x_2|+|x_4|+|x_5|)} \theta([x_1,x_2,x_4],x_5)f(x_3)\nonumber\\
      &&- (-1)^{|x_3|(|x_1|+|x_2|)}f([x_3,[x_1,x_2,x_4],x_5])\nonumber\\
       &&+\delta(-1)^{|x_3|(|x_1|+|x_2|)+|f|(|x_1|+|x_2|+|x_3|+|x_4|} D(x_3,[x_1,x_2,x_4])f(x_5)\nonumber\\
       &&-\delta(-1)^{(|f|+|x_1|+|x_2|)(|x_3|+|x_5|)+|x_4||x_5|} \theta(x_3,x_5)f([x_1,x_2,x_4]),\\
(3.10)&=&\delta(-1)^{(|x_1|+|x_2|)(|x_3|+|x_4|)}(d^{1}f)(x_3,x_4,[x_1,x_2,x_5]) \nonumber\\
      &=&\delta(-1)^{(|x_1|+|x_2|)(|x_3|+|x_4|)}
         \{(-1)^{(|f|+|x_3|)(|x_1|+|x_2|+|x_4|+|x_5|)}\theta(x_4,[x_1,x_2,x_5])f(x_3)\nonumber\\
       &&-f([x_3,x_4,[x_1,x_2,x_5]])+\delta(-1)^{|f|(|x_3|+|x_4|)}D(x_3,x_4))f([x_1,x_2,x_5]) \nonumber\\
       &&-\delta (-1)^{|x_4|(|x_1|+|x_2|+|x_5|)+|f|(|x_1|+|x_2|+|x_3|+|x_5|)}  \theta(x_3,[x_1,x_2,x_5])f(x_4) \}.
\end{eqnarray}
By (3.11)$-$(3.17), we have
\begin{eqnarray*}
&&d^{3}(d^{1}f)(x_1,x_2,x_3,x_4,x_5)\\
&=&-f([x_1,x_2,x_3],x_4,x_5)-(-1)^{|x_3|(|x_1|+|x_2|)}f(x_3,[x_1,x_2,x_4],x_5)\\
  &&-\delta(-1)^{(|x_1|+|x_2|)(|x_3|+|x_4|)}f(x_3,x_4,[x_1,x_2,x_5])+f(x_1,x_2,[x_3,x_4,x_5])\\
  &&+(-1)^{(|f|+|x_1|)(|x_2|+|x_3|+|x_4|+|x_5|)}\Lambda_{1} f(x_1)
   -\delta(-1)^{(|f|+|x_2|)(|x_3|+|x_4|+|x_5|)+|f||x_1|}\Lambda_{2} f(x_2)\\
    &&+(-1)^{|f|(|x_1|+|x_2|+|x_4|+|x_5|)+|x_3|(|x_4|+|x_5|)}\Lambda_{3} f(x_3)
   -(-1)^{|f|(|x_1|+|x_2|+|x_3|+|x_5|)+|x_4||x_5|}\Lambda_{4} f(x_4)\\
   &&+\delta(-1)^{|f|(|x_1|+|x_2|+|x_3|+|x_4|)}\Lambda_{5} f(x_5)\\
&=&(-1)^{(|f|+|x_1|)(|x_2|+|x_3|+|x_4|+|x_5|)}\Lambda_{1} f(x_1)
   -\delta(-1)^{(|f|+|x_2|)(|x_3|+|x_4|+|x_5|)+|f||x_1|}\Lambda_{2} f(x_2)\\
    &&+(-1)^{|f|(|x_1|+|x_2|+|x_4|+|x_5|)+|x_3|(|x_4|+|x_5|)}\Lambda_{3} f(x_3)
   -(-1)^{|f|(|x_1|+|x_2|+|x_3|+|x_5|)+|x_4||x_5|}\Lambda_{4} f(x_4)\\
   &&+\delta(-1)^{|f|(|x_1|+|x_2|+|x_3|+|x_4|)}\Lambda_{5} f(x_5),
\end{eqnarray*}
where
\begin{eqnarray*}
\Lambda_{1}
&=&(-1)^{(|x_2|+|x_3|)(|x_4|+|x_5|)}\theta(x_4,x_5)\theta(x_2,x_3)
   -\delta(-1)^{|x_2||x_3|+(|x_2|+|x_4|)|x_5|}\theta(x_3,x_5)\theta(x_2,x_4)\\
  &&-\theta(x_2,[x_3,x_4,x_5])+\delta(-1)^{|x_2|(|x_3|+|x_4|)}D(x_3,x_4)\theta(x_2,x_5)=0,\mbox{~by~Eq.~(3.1)}\\
\Lambda_{2}
&=& (-1)^{(|x_1|+|x_3|)(|x_4|+|x_5|)}\theta(x_4,x_5)\theta(x_1,x_3)
   -\delta(-1)^{|x_1||x_3|+|x_5|(|x_1|+|x_4|)}\theta(x_3,x_5)\theta(x_1,x_4)\\
    &&-\theta(x_1,[x_3,x_4,x_5])+\delta(-1)^{|x_1|(|x_3|+|x_4|)}D(x_3,x_4)\theta(x_1,x_5)=0,\mbox{~by~Eq.~(3.1)}\\
\Lambda_{3}
&=&  (-1)^{(|x_1|+|x_2|)(|x_4|+|x_5|)}\theta(x_4,x_5)D(x_1,x_2)
     -\delta D(x_1,x_2)  \theta(x_4,x_5)\\
     &&+\theta([x_1,x_2,x_4],x_5)
     +\delta (-1)^{|x_4|(|x_1|+|x_2|)}\theta(x_4,[x_1,x_2,x_5])=0, \mbox{~by~Eq.~(3.2)}
      \end{eqnarray*}
 \begin{eqnarray*}
\Lambda_{4}
&=& (-1)^{(|x_1|+|x_2|)(|x_3|+|x_5|)}\theta(x_3,x_5)D(x_1,x_2)
    -\delta D(x_1,x_2)  \theta(x_3,x_5)\\
     &&+\theta([x_1,x_2,x_3],x_5)
      +\delta (-1)^{|x_3|(|x_1|+|x_2|)}\theta(x_3,[x_1,x_2,x_5])=0, \mbox{~by~Eq.~(3.2)}\\
\Lambda_{5}
&=& D([x_1,x_2,x_3],x_4)+ (-1)^{|x_3|(|x_1|+|x_2|)}D(x_3,[x_1,x_2,x_4]) \\
   &&-\delta  D(x_1,x_2)D(x_3,x_4) +(-1)^{(|x_1|+|x_2|)(|x_3|+|x_4|)}D(x_3,x_4)D(x_1,x_2)=0. \mbox{~by~Eq.~(3.3)}
\end{eqnarray*}
It follows that $ d^{3}(d^{1}f)(x_1,x_2,x_3,x_4,x_5)=0$, as desired.
And this finishes the proof.$\hfill \Box$
\medskip

For $n=1,2$, the map $f\in  C_{\delta}^{n}(T,V)$ is called an $n$-cocycle if $d^nf=0$.
We denote by $Z_{\delta}^{n}(T,V)$ the subspace spanned by $n$-cocycles
and $B_{\delta}^{n}(T,V)=d^{n-2}C_{\delta}^{n-2}(T,V)$.
By Theorem 3.7, $B_{\delta}^{n}(T,V)$ is a subspace of $Z_{\delta}^{n}(T,V)$.
Therefore, we can define a cohomology  space  $H_{\delta}^{n}(T,V)$ of the $\delta$-Jordn Lie supertriple system $T$
as the factor space $Z_{\delta}^{n}(T,V)/B_{\delta}^{n}(T,V).$

\section{ 1-Parameter formal deformations of  $\delta$-Jordan Lie supertriple systems}
\def\theequation{\arabic{section}. \arabic{equation}}
\setcounter{equation} {0}

 Let $T$  be a $\delta$-Jordan Lie supertriple system and $k[[t]]$  the power series ring in one variable $t$ with coefficients in $k$.
Assume that $T[[t]]$ is the set of formal series whose coefficients are elements of the vector space $T$.
\medskip

\noindent{\bf Definition 4.1.}
Let $T$  be a $\delta$-Jordan Lie supertriple system.
 A 1-parameter formal deformations of $T$  is a formal power series
 $f_{t}: T[[t]]\times T[[t]] \times T[[t]]\rightarrow T[[t]]$  given by
\begin{eqnarray*}
f_{t}=\sum_{i\geq 0}f_{i}(x_1,x_2,x_3)t^{i}=f_{0}(x_1,x_2,x_3)+f_{1}(x_1,x_2,x_3)t+f_{2}(x_1,x_2,x_3)t^{2}+\cdots,
\end{eqnarray*}
where each $f_{i}$ is a $k$-trilinear map $f_{i}: T\times T \times T\rightarrow T$  (extended to be $k[[t]]$-trilinear)
and $f_{0}(x_1,x_2,x_3)=[x_1,x_2,x_3]$, satisfying the  following axioms:
\begin{eqnarray}
&&|f_t(x_1,x_2,x_3)|=|x_1|+|x_2|+|x_3|;\\
&&f_t(x_2,x_1,x_3)=-\delta(-1)^{|x_1||x_2|}f_t(x_1,x_2,x_3);\\
&&(-1)^{|x_1||x_3|}f_t(x_1,x_2,x_3)+(-1)^{|x_2||x_1|}f_t(x_2,x_3,x_1)+(-1)^{|x_3||x_2|}f_t(x_3, x_1,x_2)=0;~~~~~~~~\\
&&f_t(x_1,x_2,f_t(x_3,x_4,x_5))=(-1)^{|x_3|(|x_1|+|x_2|)}f_t(x_3, f_t(x_1,x_2,x_4),x_5)\nonumber\\
 &&~~~~~+f_t(f_t(x_1,x_2,x_3),x_4,x_5)+\delta(-1)^{(|x_1|+|x_2|)(|x_3|+|x_4|)}f_t(x_3,x_4, f_t(x_1,x_2,x_5)).
\end{eqnarray}

\noindent{\bf Remark 4.2.}
Equations $(4.1)-(4.4)$ are equivalent to ($n=0,1,2,\cdots$)
 \begin{eqnarray}
&&|f_i(x_2,x_1,x_3)|=|x_1|+|x_2|+|x_3|;\\
&&f_i(x_2,x_1,x_3)=-\delta(-1)^{|x_1||x_2|}f_i(x_1,x_2,x_3);\\
&&(-1)^{|x_1||x_3|}f_i(x_1,x_2,x_3)+(-1)^{|x_2||x_1|}f_i(x_2,x_3,x_1)+(-1)^{|x_3||x_2|}f_i(x_3, x_1,x_2)=0;~~~~~~~~~\\
&&\sum_{i+j=n}f_i(x_1,x_2,f_j(x_3,x_4,x_5))=\sum_{i+j=n}((-1)^{|x_3|(|x_1|+|x_2|)}f_i(x_3, f_j(x_1,x_2,x_4),x_5)\nonumber\\
  &&~~~~~~+f_i(f_j(x_1,x_2,x_3),x_4,x_5)
+\delta(-1)^{(|x_1|+|x_2|)(|x_3|+|x_4|)}f_i(x_3,x_4, f_j(x_1,x_2,x_5))).
\end{eqnarray}

 Furthermore, we can rewrite  the deformation Equation (4.8) by the equality
 $\sum_{i+j=n}f_if_j\\=0,$
where
 \begin{eqnarray*}
&&f_if_j(x_1,x_2, x_3,x_4,x_5)\\
&=&-f_i(x_1,x_2,f_j(x_3,x_4,x_5))+(-1)^{|x_3|(|x_1|+|x_2|)}f_i(x_3, f_j(x_1,x_2,x_4),x_5)\\
 &&+f_i(f_j(x_1,x_2,x_3),x_4,x_5)+\delta(-1)^{(|x_1|+|x_2|)(|x_3|+|x_4|)}f_i(x_3,x_4, f_j(x_1,x_2,x_5)).
\end{eqnarray*}

When $n=1$, Eq. (4.8) is equivalent to $f_0f_1+f_1f_0=0$.
When  $n\geq2$, Eq. (4.8) is equivalent to
$-(f_0f_n+f_nf_0)=f_1f_{n-1}+f_2f_{n-2}+\cdots+f_{n-1}f_1.$
 \medskip

 By Example 3.2, $(T,\theta)$ is the adjoint representation of $T$ itself, where
 $\theta(a,b)(x)=(-1)^{|x|(|a|+|b|)}[x,a,b]$ and $D(a,b)(x)=\delta[a,b,x]$.
 It  is easy to see that $f_i\in C_{\delta}^{3}(T,V)$  and therefore $f_if_j\in C_{\delta}^{5}(T,V)$.
 Since $f_{0}(x_1,x_2,x_3)=[x_1,x_2,x_3]$, we have
  \begin{eqnarray*}
&&f_0f_1(x_1,x_2, x_3,x_4,x_5)\\
&=&-f_0(x_1,x_2,f_1(x_3,x_4,x_5))+(-1)^{|x_3|(|x_1|+|x_2|)}f_0(x_3, f_1(x_1,x_2,x_4),x_5)\\
 &&+f_0(f_1(x_1,x_2,x_3),x_4,x_5)+\delta(-1)^{(|x_1|+|x_2|)(|x_3|+|x_4|)}f_0(x_3,x_4, f_1(x_1,x_2,x_5))\\
&=&-[x_1,x_2,f_1(x_3,x_4,x_5)]+(-1)^{|x_3|(|x_1|+|x_2|)}[x_3, f_1(x_1,x_2,x_4),x_5]\\
 &&+[f_1(x_1,x_2,x_3),x_4,x_5]+\delta(-1)^{(|x_1|+|x_2|)(|x_3|+|x_4|)}[x_3,x_4, f_1(x_1,x_2,x_5)]\\
&=&-\delta D(x_1,x_2)f(x_3,x_4,x_5)-\delta(-1)^{(|x_1|+|x_2|)(|x_3|+|x_5|)+|x_4||x_5|}D(x_3,x_4)f(x_1,x_2,x_5)\\
 &&+(-1)^{(|x_1|+|x_2|+|x_3|)(|x_4|+|x_5|)}\theta(x_4,x_5)f(x_1,x_2,x_3)\\
 &&+(-1)^{(|x_1|+|x_2|)(|x_3|+|x_4|)}D(x_3,x_4)f(x_1,x_2,x_5)
\end{eqnarray*}
 Similarly, we have
 \begin{eqnarray*}
&&f_1f_0(x_1,x_2, x_3,x_4,x_5)\\
&=&-f_1(x_1,x_2,f_0(x_3,x_4,x_5))+(-1)^{|x_3|(|x_1|+|x_2|)}f_1(x_3, f_0(x_1,x_2,x_4),x_5)\\
 &&+f_1(f_0(x_1,x_2,x_3),x_4,x_5)+\delta(-1)^{(|x_1|+|x_2|)(|x_3|+|x_4|)}f_1(x_3,x_4, f_0(x_1,x_2,x_5))\\
&=&-f_1(x_1,x_2,[x_3,x_4,x_5])+(-1)^{|x_3|(|x_1|+|x_2|)}f_1(x_3, [x_1,x_2,x_4],x_5)\\
 &&+f_1([x_1,x_2,x_3],x_4,x_5)+\delta(-1)^{(|x_1|+|x_2|)(|x_3|+|x_4|)}f_1(x_3,x_4, [x_1,x_2,x_5]).
\end{eqnarray*}
It follows that
  \begin{eqnarray*}
(f_0f_1+f_1f_0)(x_1,x_2, x_3,x_4,x_5)
=d^{3}f_{1}(x_1,x_2, x_3,x_4,x_5).
\end{eqnarray*}
 Therefore, we deduce  $d^{3}f_{1}=0$ since $f_0f_1+f_1f_0=0$.
 Also we can obtain $-d^3f_n=f_1f_{n-1}+f_2f_{n-2}+\cdots+f_{n-1}f_1.$
And $f_{1}$ is called the infinitesimal  deformation of $f_{t}$.
\medskip

\noindent{\bf Definition 4.3.}
Let $T$  be a $\delta$-Jordan Lie supertriple system.
 Two 1-parameter formal deformations $f_{t}$ and $f'_{t}$  of $T$ are said to be equivalent,  denoted by $f_{t}\sim f'_{t}$,
 if there exists a formal isomorphism of  $k[[t]]$-modules
  \begin{eqnarray*}
\phi_{t}(x)=\sum_{i\geq 0}\phi_{i}(x)t^{i}:(T[[t]],f_{t},\delta)\rightarrow (T[[t]],f'_{t},\delta),
\end{eqnarray*}
where $\phi_{0}=id_T$, $\phi_{i}:T \rightarrow T$ is an $k-$linear map (extended to be $k[[t]]-$linear) such that
 \begin{eqnarray*}
\phi_{t}f_t(x_1,x_2, x_3)=f'_t(\phi_{t}(x_1),\phi_{t}(x_2), \phi_{t}(x_3)), ~\forall x_1,x_2, x_3\in T.
\end{eqnarray*}

In particular, if $f_1=f_2=\cdots =0,$ then  $f_1=f_0$ is called the null deformation.
If  $f_t\sim f_0$, then $f_t$ is called the trivial deformation.
If every 1-parameter formal deformation $f_t$ is trivial, then $T$ is called an analytically rigid $\delta$-Jordan Lie supertriple system.

\medskip

\noindent{\bf Theorem 4.4.}
Let $f_{t}=\sum_{i\geq 0}f_{i}(x_1,x_2,x_3)t^{i}$ and $f'_{t}=\sum_{i\geq 0}f'_{i}(x_1,x_2,x_3)t^{i}$ be
 two equivalent 1-parameter formal deformations  of $T$.
Then the infinitesimal deformations $f_1$ and $f'_1$
belong to the same cohomology class in  $H_{\delta}^{3}(T,T).$

\medskip

\noindent{\bf Proof.}
%It is enough to prove that $f_1-f'_1\in B_{\delta}^{3}(T,T)$.
By the assumption that $f_1$ and $f'_1$ are equivalent, there exists  a formal isomorphism $\phi_{t}(x)=\sum_{i\geq 0}\phi_{i}(x)t^{i}$
 of  $k[[t]]$-modules satisfying
  \begin{eqnarray*}
\sum_{i\geq 0}\phi_i(\sum_{j\geq 0}f_j(x_1,x_2, x_3)t^{j})
=\sum_{i\geq 0}f'_i(\sum_{k\geq 0}\phi_{k}(x_1)t^{k},\sum_{l\geq 0}\phi_{l}(x_2)t^{l}, \sum_{m\geq 0}\phi_{m}(x_3)t^{m})t^{i},
\end{eqnarray*}
for any $x_1,x_2, x_3\in T.$
Comparing with the coefficients of $t^1$ for two sides of the above equation, we have
    \begin{eqnarray*}
&&f_1(x_1,x_2, x_3)+\phi_{1}([x_1,x_2, x_3])\\
&=&f'_1(x_1,x_2, x_3)+[\phi_{1}(x_1),x_2, x_3]+[x_1,\phi_{1}(x_2), x_3]+[x_1,x_2, \phi_{1}(x_3)].
\end{eqnarray*}
It follows that
  \begin{eqnarray*}
&&f_1(x_1,x_2, x_3)-f'_1(x_1,x_2, x_3)\\
&=&[\phi_{1}(x_1),x_2, x_3]+[x_1,\phi_{1}(x_2), x_3]+[x_1,x_2, \phi_{1}(x_3)]-\phi_{1}([x_1,x_2, x_3])\\
&=&(-1)^{|x_1|(|x_2|+|x_3|)}\theta(x_2,x_3)\phi_{1}(x_1)-\phi_{1}([x_1,x_2,x_3])\\
    &&+\delta D(x_1,x_2)\phi_{1}(x_3)
    -\delta(-1)^{|x_2||x_3|}\theta(x_1,x_3)\phi_{1}(x_2)\\
&=& d^{1}\phi_{1} (x_1,x_2, x_3).
\end{eqnarray*}
So $f_1-f'_1=d^{1}\phi_{1}\in B_{\delta}^{3}(T,T),$  as dsired.
The proof is completed.$\hfill \Box$
\medskip

\noindent{\bf Theorem 4.5.}
Let $T$  be a $\delta$-Jordan Lie supertriple system
with $H_{\delta}^{3}(T,T)=0$, then $T$ is   analytically rigid.
\medskip

\noindent{\bf Proof.}
Let $f_{t}=\sum_{i\geq 0}f_{i}t^{i}$ be a 1-parameter formal deformation  of $T$.
Then
 $d^3f_n=f_1f_{n-1}+f_2f_{n-2}+\cdots+f_{n-1}f_1=0.$
 By the assumption $H_{\delta}^{3}(T,T)=0$, we have  $f_n\in Z_{\delta}^{3}(T,T)=B_{\delta}^{3}(T,T)$,
that is, there exits $g_n\in C_{\delta}^{1}(T,T)$ such that $f_n=d^{1}g^n$.

Set $\phi_t=id_T-g_nt^n,$ then
\begin{eqnarray*}
&&\phi_t(id_T+g_nt^n+g_n^2t^{2n}+g_n^3t^{3n}+\cdots)\\
&&~~~~~=(id_T-g_nt^n)(id_T+g_nt^n+g_n^2t^{2n}+g_n^3t^{3n}+\cdots)\\
&&~~~~~=(id_T+g_nt^n+g_n^2t^{2n}+g_n^3t^{3n}+\cdots)-(g_nt^n+g_n^2t^{2n}+g_n^3t^{3n}+\cdots)\\
&&~~~~~=id_T.
\end{eqnarray*}
Similarly, one may check that $(id_T+g_nt^n+g_n^2t^{2n}+g_n^3t^{3n}+\cdots)\phi_t=id_T.$
So $\phi_t:(T[[t]],f_{t},\delta)\rightarrow (T[[t]],f'_{t},\delta)$ is a linear isomorphism.
Thus we can define another  1-parameter formal deformation  by $\phi_t^{-1}$ in the form of
$$f'_{t}(x_1,x_2,x_3)=\phi_t^{-1}f_{i}(\phi_t(x_1),\phi_t(x_2),\phi_t(x_3)).$$
Obviously, $f_{t}\sim f'_{t}$.
Set  $f'_{t}=\sum_{i\geq 0}f'_{i}t^{i}$, then we have
\begin{eqnarray*}
&&\sum_{i\geq 0}f'_{i}(x_1,x_2,x_3)t^{i}-\sum_{i\geq 0}g_nf_{i}(x_1,x_2,x_3)t^{i+n}\\
&=&f_0(x_1,x_2,x_3)
   -\{f_0(g_n(x_1),x_2,x_3)+f_0(x_1,g_n(x_2),x_3)+f_0(x_1,x_2,g_n(x_3))\}t^n\\
   &&+\{f_0(g_n(x_1),g_n(x_2),x_3)+f_0(x_1,g_n(x_2),g_n(x_3))+f_0(x_1,g_n(x_2),g_n(x_3))\}t^{2n}\\
   &&+f_0(g_n(x_1),g_n(x_2),g_n(x_3))t^{3n}+\sum_{i\geq n}f_{i}(x_1,x_2,x_3)t^{i}\\
   &&-\sum_{i\geq n}\{f_i(g_n(x_1),x_2,x_3)+f_i(x_1,g_n(x_2),x_3)+f_i(x_1,x_2,g_n(x_3))\}t^{i+n}\\
   &&+\sum_{i\geq n}\{f_i(g_n(x_1),g_n(x_2),x_3)+f_i(x_1,g_n(x_2),g_n(x_3))+f_i(g_n(x_1),x_2,g_n(x_3))\}t^{i+2n}\\
    &&-\sum_{i\geq n}f_i(g_n(x_1),g_n(x_2),g_n(x_3))t^{i+3n}.
\end{eqnarray*}
By the above equation, it follows that
\begin{eqnarray*}
&&f'_0(x_1,x_2,x_3)=f_0(x_1,x_2,x_3)=[x_1,x_2,x_3],\\
&&f'_1(x_1,x_2,x_3)=f'_1(x_1,x_2,x_3)=\cdots=f'_{n-1}(x_1,x_2,x_3)=0,\\
&&f'_n(x_1,x_2,x_3)-g_n([x_1,x_2,x_3])\\
&&~~~~~~~=f_n(x_1,x_2,x_3)-[g_n(x_1),x_2,x_3]-[x_1,g_n(x_2),x_3]-[x_1,x_2,g_n(x_3)]\\
&&~~~~~~~=f_n(x_1,x_2,x_3)-(-1)^{|x_1|(|x_2|+|x_3|)}\theta(x_2,x_3)g_n(x_1)\\
&&~~~~~~~~~~~             +\delta(-1)^{|x_2||x_3|}\theta(x_1,x_3)g_n(x_2)-\delta D(x_1,x_2)g_n(x_3).
\end{eqnarray*}
Therefore, we deduce
\begin{eqnarray*}
f'_n(x_1,x_2,x_3)=f_n(x_1,x_2,x_3)-d^1g_n(x_1,x_2,x_3)=0.
\end{eqnarray*}
It follows that $f'_t=f_0+\sum_{i\geq n+1}f'_it^i$.
 By induction, we have $f_t\sim f_0$, that is, $T$ is   analytically rigid.
 The proof is finished.
$\hfill \Box$

\section{Nijenhuis operators of $\delta$-Jordan Lie supertriple systems}
\def\theequation{\arabic{section}. \arabic{equation}}
\setcounter{equation} {0}

In this section, we introduce the notion of
Nijenhuis operators for $\delta$-Jordan Lie supertriple systems.
Also, we give trivial deformations of this kind of operators.
\medskip

Let $T$  be a $\delta$-Jordan Lie supertriple system and $\psi: T\times T \times T\rightarrow T$  be an even trilinear map.
Consider a $\lambda$-parametrized family of linear operations:
\begin{eqnarray*}
[x_{1},x_{2},x_{3}]_{\lambda}=[x_{1},x_{2},x_{3}]+\lambda \psi(x_{1},x_{2},x_{3}),
\end{eqnarray*}
for any $x_{1},x_{2},x_{3}\in T$,  where $\lambda$ is a formal variable.
\medskip

If $[\c,\c,\c]_{\lambda}$ endow $T$ with a $\delta$-Jordan Lie supertriple system structure which is denoted by $T_\lambda$ , then we call that $\psi$ generates a $\lambda$-parameter infinitesimal deformation of the $\delta$-Jordan Lie supertriple system  $T$.
\medskip

\noindent{\bf Theorem 5.1.}
Let $T$  be a $\delta$-Jordan Lie supertriple system.
Then  $T_\lambda$ is a $\delta$-Jordan Lie supertriple system if and only if

 (i) $\psi$  itself defines a $\delta$-Jordan Lie supertriple system structure on $T$;

 (ii) $\psi$ is a 3-cocycle of $T$.
\medskip

\noindent{\bf Proof.}
Assume that $T_{\lambda}$ is a $\delta$-Jordan Lie supertriple system.
For any $x_{1},x_{2},x_{3}\in T$, we have
\begin{eqnarray*}
&&[x_{2},x_{1},x_{3}]_{\lambda}
=[x_{2},x_{1},x_{3}]+\lambda\psi(x_{2},x_{1},x_{3})
=-\delta(-1)^{|x_{1}||x_{2}|}[x_{1},x_{2},x_{3}]+\lambda\psi(x_{2},x_{1},x_{3}),\\
&&[x_{2},x_{1},x_{3}]_{\lambda}
=-\delta(-1)^{|x_{1}||x_{2}|}[x_{1},x_{2},x_{3}]_{\lambda}
=-\delta(-1)^{|x_{1}||x_{2}|}[x_{1},x_{2},x_{3}]-\delta(-1)^{|x_{1}||x_{2}|}\lambda \psi(x_{1},x_{2},x_{3}).
\end{eqnarray*}
It follows that
\begin{eqnarray}
\psi(x_{2},x_{1},x_{3})=-\delta(-1)^{|x_{1}||x_{2}|}\psi(x_{1},x_{2},x_{3}).
\end{eqnarray}

For Eq. (2.3), we have
\begin{eqnarray}
&&\lambda((-1)^{|x_{1}||x_{3}|}\psi(x_{1},x_{2},x_{3})+(-1)^{|x_{2}||x_{1}|}\psi(x_{2},x_{3},x_{1})+(-1)^{|x_{3}||x_{2}|}\psi(x_{3},x_{1},x_{2}))\nonumber\\
&=&(-1)^{|x_{1}||x_{3}|}[x_{1},x_{2},x_{3}]_{\lambda}+(-1)^{|x_{2}||x_{1}|}[x_{2},x_{3},x_{1}]_{\lambda}+(-1)^{|x_{3}||x_{2}|}[x_{3},x_{1},x_{2}]_{\lambda}\nonumber\\
&&-(-1)^{|x_{1}||x_{3}|}[x_{1},x_{2},x_{3}]-(-1)^{|x_{2}||x_{1}|}[x_{2},x_{3},x_{1}]-(-1)^{|x_{3}||x_{2}|}[x_{3},x_{1},x_{2}]\nonumber\\
&=&0,
\end{eqnarray}
as desired. The last equality holds since  $T$ and  $T_\lambda$  are both $\delta$-Jordan Lie supertriple systems.

For Eq. (2.3), we  take  $x_{1},x_{2},x_{3},x_{4},x_{5}\in T$ and calculate
\begin{eqnarray*}
&&[x_{1},x_{2},[x_{3},x_{4},x_{5}]_{\lambda}]_{\lambda}\\
&=&[x_{1},x_{2},[x_{3},x_{4},x_{5}]+\lambda\psi(x_{3},x_{4},x_{5})]_{\lambda}
\end{eqnarray*}
\begin{eqnarray*}
&=&[x_{1},x_{2},[x_{3},x_{4},x_{5}]+\lambda\psi(x_{3},x_{4},x_{5})]+\lambda\psi(x_{1},x_{2},[x_{3},x_{4},e]+\lambda\psi(x_{3},x_{4},x_{5}))\\
&=&[x_{1},x_{2},[x_{3},x_{4},x_{5}]]+\lambda([x_{1},x_{2},\psi(x_{3},x_{4},x_{5})]
   +\psi(x_{1},x_{2},[x_{3},x_{4},x_{5}]))\\
   &&+\lambda^{2}\psi(x_{1},x_{2},\psi(x_{3},x_{4},x_{5})).
\end{eqnarray*}
By similar calculation, we  have
\begin{eqnarray*}
&&[[x_{1},x_{2},x_{3}]_{\lambda},x_{4},x_{5}]_{\lambda}
=[[x_{1},x_{2},x_{3}],x_{4},x_{5}]+\lambda([\psi(x_{1},x_{2},x_{3}),x_{4},x_{5}]\\
     &&~~~~~~~~~~~~~~~~~~~~~~~~~~~~~~+\psi([x_{1},x_{2},x_{3}],x_{4},x_{5}))+\lambda^{2}\psi(\psi(x_{1},x_{2},x_{3}),x_{4},x_{5}),\\
&&[x_{3},[x_{1},x_{2},x_{4}]_{\lambda},x_{5}]_{\lambda}
=[x_{3},[x_{1},x_{2},x_{4}],x_{5}]+\lambda([x_{3},\psi(x_{1},x_{2},x_{4}),x_{5}]\\
     &&~~~~~~~~~~~~~~~~~~~~~~~~~~~~~~+\psi(x_{3},[x_{1},x_{2},x_{4}],x_{5}))+\lambda^{2}\psi(x_{3},\psi(x_{1},x_{2},x_{4}),x_{5}),\\
&&[x_{3},x_{4},[x_{1},x_{2},x_{5}]_{\lambda}]_{\lambda}
=[x_{3},x_{4},[x_{1},x_{2},x_{5}]]+\lambda([x_{3},x_{4},\psi(x_{1},x_{2},x_{5})]\\
     &&~~~~~~~~~~~~~~~~~~~~~~~~~~~~~~+\psi(x_{3},x_{4},[x_{1},x_{2},x_{5}]))+\lambda^{2}\psi(x_{3},x_{4},\psi(x_{1},x_{2},x_{5})).
\end{eqnarray*}
It follows that
\begin{eqnarray*}
&&[[x_{1},x_{2},x_{3}]_{\lambda},x_{4},x_{5}]_{\lambda}+(-1)^{|x_{3}|(|x_{1}|+|x_{2}|)}[x_{3},[x_{1},x_{2},x_{4}]_{\lambda},x_{5}]_{\lambda}\\
      &&+\delta(-1)^{(|x_{1}|+|x_{2}|)(|x_{3}|+|x_{4}|)}[x_{3},x_{4},[x_{1},x_{2},x_{5}]_{\lambda}]_{\lambda}\\
&=&[[x_{1},x_{2},x_{3}],x_{4},e]+(-1)^{|x_{3}|(|x_{1}|+|x_{2}|)}[x_{3},[x_{1},x_{2},x_{4}],x_{5}]\\
      &&+\delta(-1)^{(|x_{1}|+|x_{2}|)(|x_{3}|+|x_{4}|)}[x_{3},x_{4},[x_{1},x_{2},x_{5}]]\\
      &&+\lambda\{[\psi(x_{1},x_{2},x_{3}),x_{4},e]+\psi([x_{1},x_{2},x_{3}],x_{4},x_{5})
      +(-1)^{|x_{3}|(|x_{1}|+|x_{2}|)}[x_{3},\psi(x_{1},x_{2},x_{4}),x_{5}]\\
      &&+(-1)^{|x_{3}|(|x_{1}|+|x_{2}|)}\psi(x_{3},[x_{1},x_{2},x_{4}],x_{5})
      +\delta(-1)^{(|x_{1}|+|x_{2}|)(|x_{3}|+|x_{4}|)}[x_{3},x_{4},\psi(x_{1},x_{2},x_{5})]\\
      &&+\delta(-1)^{(|x_{1}|+|x_{2}|)(|x_{3}|+|x_{4}d|)}\psi(x_{3},x_{4},[x_{1},x_{2},x_{5}])\}
      +\lambda^{2}\{\psi(\psi(x_{1},x_{2},x_{3}),x_{4},x_{5})\\
      &&+(-1)^{|x_{3}|(|x_{1}|+|x_{2}|)}\psi(x_{3},\psi(x_{1},x_{2},x_{4}),x_{5})
         +\delta(-1)^{(|x_{1}|+|x_{2}|)(|x_{3}|+|x_{4}|)}\psi(x_{3},x_{4},\psi(x_{1},x_{2},x_{5}))\}.
\end{eqnarray*}
Therefore, we have
\begin{eqnarray}
&&[x_{1},x_{2},\psi(x_{3},x_{4},x_{5})]+\psi(x_{1},x_{2},[x_{3},x_{4},x_{5}])\nonumber\\
&=&[\psi(x_{1},x_{2},x_{3}),x_{4},x_{5}]+\psi([x_{1},x_{2},x_{3}],x_{4},x_{5})
      +(-1)^{|x_{3}|(|x_{1}|+|x_{2}|)}[x_{3},\psi(x_{1},x_{2},x_{4}),x_{5}]\nonumber\\
      &&+(-1)^{|x_{3}|(|x_{1}|+|x_{2}|)}\psi(x_{3},[x_{1},x_{2},x_{4}],x_{5})
      +\delta(-1)^{(|x_{1}|+|x_{2}|)(|x_{3}|+|x_{4}|)}[x_{3},x_{4},\psi(x_{1},x_{2},x_{5})]\nonumber\\
      &&+\delta(-1)^{(|x_{1}|+|x_{2}|)(|x_{3}|+|x_{4}|)}\psi(x_{3},x_{4},[x_{1},x_{2},x_{5}]),\\
&&\psi(x_{1},x_{2},\psi(x_{3},x_{4},x_{5}))\nonumber\\
 &=& \psi(\psi(x_{1},x_{2},x_{3}),x_{4},x_{5})
      +(-1)^{|x_{3}|(|x_{1}|+|x_{2}|)}\psi(x_{3},\psi(x_{1},x_{2},x_{4}),x_{5})\nonumber\\
      &&+\delta(-1)^{(|x_{1}|+|x_{2}|)(|x_{3}|+|x_{4}|)}\psi(x_{3},x_{4},\psi(x_{1},x_{2},x_{5})).
\end{eqnarray}
By Eq. (5.4),  $\psi$ satisfies  Eq. (2.3). So $\psi$ defines a $\delta$-Jordan Lie supertriple system structure on $T$.

Since $\theta(x_{1},x_{2})(x)=(-1)^{|x|(|x_{1}|+|x_{2}|)}[x,x_{1},x_{2}]$ and $D(x_{1},x_{2})(x)=\delta[x_{1},x_{2},x]$,
 we can rewrite  Eq. (5.2) as follows:
\begin{eqnarray}
0&=&-\delta D(x_{1},x_{2})\psi(x_{3},x_{4},x_{5})-\psi(x_{1},x_{2},[x_{3},x_{4},x_{5}])\nonumber\\
&&+(-1)^{(|x_{1}|+|x_{2}|+|x_{3}|)(|x_{4}|+|x_{5}|)}\theta(x_{4},x_{5})\psi(x_{1},x_{2},x_{3})+\psi([x_{1},x_{2},x_{3}],x_{4},x_{5})\nonumber\\
     &&-\delta(-1)^{|x_{3}|(|x_{1}|+|x_{2}|)}\theta(x_{3},e)\psi(x_{1},x_{2},x_{4})
       +(-1)^{|x_{3}|(|x_{1}|+|x_{2}|)}\psi(x_{3},[x_{1},x_{2},x_{4}],x_{5})\nonumber\\
      &&+\delta(-1)^{(|x_{1}|+|x_{2}|)(|x_{3}|+|x_{4}|)}D(x_{3},x_{4})\psi(x_{1},x_{2},x_{5})\nonumber\\
      &&+\delta(-1)^{(|x_{1}|+|x_{2}|)(|x_{3}|+|x_{4}|)}\psi(x_{3},x_{4},[x_{1},x_{2},x_{5}])\nonumber\\
&=&d^{3}\psi(x_{1},x_{2},x_{3},x_{4},x_{5}).
\end{eqnarray}
The last equality holds since  $\psi$  is an even trilinear map.
So $d^{3}\psi=0$, as required.

Conversely, if  $\psi$ satisfies conditions (i) and (ii), it is easy to see that
$T_\lambda$ is a $\delta$-Jordan Lie supertriple system from Eqs. $(5.1)-(5.4).$$\hfill \Box$
\medskip

\noindent{\bf Definition 5.2.}
A deformation is said to be trivial if there exists a linear map $N:T\rightarrow T$
such that for all $\lambda$, $\varphi_{\lambda}=id+\lambda N:T_\lambda\rightarrow T$
satisfies
\begin{eqnarray}
\varphi_{\lambda}[x_1,x_2,x_3]_\lambda=[\varphi_{\lambda}x_1,\varphi_{\lambda}x_2,\varphi_{\lambda}x_3],
\end{eqnarray}
for any $x_1,x_2,x_3\in T$.
\medskip

The left hand side of Eq. (5.6) equals to
\begin{eqnarray*}
[x_1,x_2,x_3]+\lambda\{\psi(x_1,x_2,x_3)+N[x_1,x_2,x_3]\}+\lambda^2 N\psi(x_1,x_2,x_3).
\end{eqnarray*}

The right hand side of Eq. (5.6) equals to
\begin{eqnarray*}
&&[x_1,x_2,x_3]+
\lambda\{[N x_1,x_2,x_3]+[x_1,N x_2,x_3]+[x_1,x_2,N x_3]\}\\
&&~~+\lambda^2 \{[N x_1,N x_2,x_3]+[x_1,N x_2,N x_3]+[N x_1,x_2,N x_3)]\}
+\lambda^3[N x_1,N x_2,N x_3].
\end{eqnarray*}

Therefore, by Eq. (5.6), we have
\begin{eqnarray}
0&=&[N x_1,N x_2,N x_3],\\
N\psi(x_1,x_2,x_3)&=&[N x_1,N x_2,x_3]+[x_1,N x_2,N x_3]+[N x_1,x_2,N x_3)],\\
\psi(x_1,x_2,x_3)&=&[N x_1,x_2,x_3]+[x_1,N x_2,x_3]+[x_1,x_2,N x_3]-N[x_1,x_2,x_3]\nonumber\\
                 &=&(-1)^{|x_1|(|x_2|+|x_3|)}\theta(x_2,x_3)N x_1-\delta(-1)^{|x_2||x_3|}\theta(x_1,x_3)N x_2\nonumber\\
                 &&~~+\delta D(x_1,x_2)N x_3-N[x_1,x_2,x_3]\nonumber\\
                 &=&d^1N(x_1,x_2,x_3).
\end{eqnarray}

By Eq. (5.8) and Eq. (5.9), we can deduce that
\begin{eqnarray}
N^2[x_1,x_2,x_3]
&=&N[N x_1,x_2,x_3]+N[x_1,N x_2,x_3]+N[x_1,x_2,N x_3]\nonumber\\
&&~~-[N x_1,N x_2,x_3]-[x_1,N x_2,N x_3]-[N x_1,x_2,N x_3)].
\end{eqnarray}

\noindent{\bf Definition 5.3.}
A linear operator $N:T\rightarrow T$ is called a Nijenhuis operator if and only if
Eq. (5.7) and Eq. (5.10) hold.
\medskip

\noindent{\bf Theorem 5.4.}
Let $N$ be a Nijenhuis operator for $T$. Then, a deformation of $T$ can be obtained by putting
\begin{eqnarray*}
\psi(x_1,x_2,x_3) &=&(-1)^{|x_1|(|x_2|+|x_3|)}\theta(x_2,x_3)N x_1-\delta(-1)^{|x_2||x_3|}\theta(x_1,x_3)N x_2\\
                 &&+\delta D(x_1,x_2)N x_3-N[x_1,x_2,x_3].
\end{eqnarray*}
Moreover, this deformation is trivial.
\medskip

\noindent{\bf Proof.}
 The proof is similar to one in the setting of  $\delta$-Jordan Lie triple system  in \cite{chen2017}.
$\hfill \Box$

\begin{center}
 {\bf ACKNOWLEDGEMENT}
 \end{center}

  The paper is  supported by
   the Anhui Provincial Natural Science Foundation (No. 1808085MA14),
   the NSF of China (No. 11761017) and the Youth Project for Natural Science Foundation of Guizhou provincial department of education (No. KY[2018]155).

\renewcommand{\refname}{REFERENCES}

\end{document}